\documentclass{article}
\usepackage{graphicx}
 \usepackage{mathptmx}
\usepackage{amsmath, amstext,amssymb,amsfonts}
\usepackage[english]{babel}
\usepackage{delarray}
\usepackage{mathptmx}
\usepackage{amsmath, amstext,amssymb,amsfonts}
\usepackage[english]{babel}
\usepackage{delarray}

\DeclareFontFamily{U}{mathx}{\hyphenchar\font45}
\DeclareFontShape{U}{mathx}{m}{n}{
      <5> <6> <7> <8> <9> <10>
      <10.95> <12> <14.4> <17.28> <20.74> <24.88>
      mathx10
      }{}
\DeclareSymbolFont{mathx}{U}{mathx}{m}{n}
\DeclareFontSubstitution{U}{mathx}{m}{n}
\DeclareMathAccent{\widecheck}{0}{mathx}{"71}
\DeclareMathAccent{\wideparen}{0}{mathx}{"75}

\numberwithin{equation}{section}

\newtheorem{thm}{Theorem}[section]
 
 \newtheorem{lem}[thm]{Lemma}
 \newtheorem{prop}[thm]{Proposition}
      \newcommand{\us}{ U_{\sigma}}     \newcommand{\bev}{B_{\qu}^{1}}
                    \newcommand{\su}{{\text{\,supp}}\,}         \newcommand{\ben}{B^1_{Q^n_{\sigma}}}         \newcommand{\bep}{B^1_{Q^n_{2\pi}} }                                 \newcommand{\qu}{Q^n_{\sigma}}        \newcommand{\sn}{\Delta_{\,n}}  \newcommand{\Bn}{B(\Rn)}
  \newcommand{\R}{\mathbb R}    \newcommand{\Co}{\mathbb C}             \newcommand{\Z}{\mathbb Z}  \newcommand{\Rn}{\mathbb{R}^n} \newcommand{\Zn}{\mathbb{Z}^n}   \newcommand{\Mn}{M({\mathbb R}^n)}  \newcommand{\Con}{\mathbb {C}^n}

\date{}
\title{ On the deterministic   property for characteristic  functions of several variables}
\author{\large{ Saulius  Norvidas}}
\date{\footnotesize Institute of Data Science and Digital Technologies, Vilnius University, \\ Akademijos str. 4, Vilnius LT-04812, Lithuania\\
 ({\rm{e-mail: norvidas{@}gmail.com}})}
\begin{document}

\maketitle
 {{ {\bf Abstract}}} \ Assume that $f$ is the characteristic function of a probability measure $\mu_f$ on $\Rn$. Let $\sigma>0$.  We study the following extrapolation problem:  under what conditions on the neighborhood of infinity $V_{\sigma}=\{x\in\Rn: |x_k|>\sigma, \  k=1,\dots, n\}$ in $\Rn$ does there exist a characteristic function $g$ on $\Rn$ such that $g=f$ on $V_{\sigma}$, but  $g\not\equiv f$? Let $\mu_f$ have a nonzero absolutely continuous part  with continuous  density $\varphi$.  In this paper  certain sufficient conditions on $\varphi$ and $V_{\sigma}$ are given  under which the latter question has an  affirmative answer. We also address the optimality of these conditions. Our results indicate that not only does the size of both  $V_{\sigma}$ and the support ${{\text{\,supp}}\,}\varphi$ matter, but also certain arithmetic properties of ${{\text{\,supp}}\,}\varphi$.

{\bf Keywords}: Characteristic function;   density function;   entire function;  probability measure;  Bernstein space

{\bf  Mathematics Subject Classification}:   32D15 - 42A38 - 60E10

\section{ Introduction }

{\large{
A function $f : \Rn\to\Co$  is said to be positive definite if the inequality
\[
\sum_{j, k=1}^{n} f(x_j-x_k)c_j\overline{c}_k\ge 0
\]
holds for all finite sets of complex numbers $c_1,\dots, c_n$ and points $x_1, \dots, x_n \in\Rn$. Any such a function $f$ satisfies
\begin{equation}
f(-x)=\overline{f(x)}\quad{\text{and}}\quad |f(x)|\le f(0)
\end{equation}
for all $x\in\Rn$. The Bochner theorem gives the description of continuous positive definite functions in terms of the Fourier transform. For this reason, let us  recall certain notions.

Let  $\Mn$  be  the  Banach algebra of bounded regular complex-valued  Borel measures $\mu$ on $\Rn$. As usual, $\Mn$ is  equipped with the total variation norm $\|\mu\|$.  For $\mu\in\Mn$,  we define the  Fourier  transform   by
\[
\hat{\mu}(x)=\int_{\Rn}e^{-i( x, t)} d\mu(t),
\]
 $ x\in\Rn$, where $( x, t)=\sum_{k=1}^{n}  x_kt_k$ is the  scalar product on $\Rn$.   We identify $L^1(\Rn)$ with the closed ideal in $\Mn$ of all measures which are absolutely continuous  with respect to the Lebesgue  measure $dt=dt_1\cdot {\dots} \cdot dt_n$ on $\Rn$.

Let   $\mu\in \Mn$ be a positive measure. If   $\|\mu\|=1$,  then in the language of probability theory,  $\mu$ and  $f(x):= \hat{\mu}(x)$, $x\in\Rn$,  are called a  probability measure and its characteristic function, respectively.  In particular, if $\mu=\varphi\, dt$ with  $\varphi\in L^1(\Rn)$   such  that  $\|\varphi\|_{L^1(\Rn)}=1$ and $\varphi\ge 0$ on $\Rn$,   then $\varphi$  is called the probability density function of $\mu$, or the probability density for short. The Bochner theorem (see, e.g., \cite[p. 57]{7}) states that a function $f:\Rn\to\Co$ is the characteristic function of  a  probability measure if and only if $f$ is continuous and positive definite on $\Rn$ with $f(0)=1$.

Given a characteristic function $f$, we are interested in the   question: Is it true that for any open subset  $U$ of  $\Rn$, $0\not\in U$,
there exists a characteristic function $g$ such that $g=f$   on $U$  but $g \not\equiv f$? Our interest to this question is initiated by a similar  problem posed by N.G. Ushakov in  (\cite[p. 276]{7}):  Is it true that for any  interval $[a, b]\subset \mathbb{R} $, $0\notin [a,b]$,  there exists the characteristic function $f$  such that  $f \not\equiv e^{-x^2/2}$, but  $ f (x) = e^{-x^2/2}$ for all $x\in [a, b]$?  A positive answer to this question  follows immediately from the following result of  Gneiting (see \cite[p.360]{1}):

{\it  Let $f : \R\to\Co$  be the characteristic function of a  distribution with a continuous and strictly positive density.  Then there exists, for each $a>0$,   a characteristic function $g$   such that $f(x)=g(x)$ if $x=0$ or $|x|\ge a$ and  $f(x)\neq g(x)$ otherwise.}

Next, in \cite[p.  361]{1}, the author says that  a characteristic function $f$ has the substitution property if, for any $a>0$, there exists a characteristic function $g$ such that $f(x)=g(x)$ for $|x|> a$ but   $f\not\equiv g$. Also in \cite[p.  361]{1},  we  can find such a conjecture:

{\it One might conjecture that any characteristic function with an absolutely continuous component has the substitution property}.

This assumption is not confirmed. Indeed, it follows from  the following our result (see \cite[p. 238]{3}):
\begin{prop}.\label{t1}
 Let $f : \R\to\Co$  be the characteristic function of the triangular (Simpson) distribution with the density function ƒ$\varphi$ defined by
\[
\varphi(t)  =\begin{array}\{{rl}.
1- |t|, &  |t|\le 1, \\
0, & \mbox{otherwise.}
\end{array}
\]
Suppose  $g : \R\to\Co$ is a characteristic function such that $f(x)=g (x)$ for all $|x|>\varrho$ with certain $\varrho> 0$. If $\varrho\le \pi/2$, then $f\equiv  g$.
\end{prop}

In this paper,  a problem of uniqueness  for an extrapolation of characteristic functions of several variables is studied.  More precisely, given  a characteristic function $f:\Rn\to\Co$, we consider characteristic extrapolations of $f$,   in a manner indicated by the above mentioned Ushakov's problem,  from    neighborhoods $U$ of  infinity to the whole $\Rn$. Any characteristic function $f:\Rn\to\Co$ satisfies  $f(-x)=\overline{f(x)}$, $ x\in\Rn$. This implies that it is enough to study the extensions of $f$ only from symmetric neighborhoods of infinity.

 Let us recall certain   notations.  Any characteristic function $f:\Rn\to\Co$ satisfies  $f(-x)=\overline{f(x)}$, $ x\in\Rn$. This implies that it is enough to study the extensions of $f$ only from symmetric neighborhoods of infinity.  For $\sigma>0$,  set $ \qu=\{x\in\Rn: |x_k|\le \sigma, k=1,\dots,n\}$. Then
 \[
V_{\sigma}=\Rn\setminus  \qu
\]
  denotes such a neighborhood.   We say that  $f$ has the $\sigma$-deterministic property if there exists no other  characteristic function $g$ such that $f(x)=g(x)$ for all $x\in V_{\sigma}$.

  Given a continuous function $\varphi:\Rn\to\Co$, we define its support $\su \varphi $  as usual. Set
  \[
  S(\varphi)=\{t\in\Rn: \ |\varphi(t)|>0\}.
   \]
   We call $ S(\varphi)$ the essential support of $\varphi$. Then $N(\varphi):=\Rn\setminus S(\varphi)$ is the zero set of $\varphi$. If $A$ and $B$ are subsets of $\Rn$ and $\omega\in\Rn$, then $A+B$ and $\omega\cdot A$ denote the sets $\{a+b:\  a\in A, b\in B\}$ and $\{(\omega_1a_1,\dots,\omega_na_n):\  a\in A\}$, respectively.  In particular, if $\omega_k\neq 0$ for all $k$, then  $(1/\omega) A=\{\bigl((1/\omega_1)a_1,\dots,(1/\omega_n)a_n\bigr):\  a=(a_1,\dots, a_n)\in A\}$. Finally, let   $\Zn$  be the usual  integer $n$-lattice, where  $\Z$ is the group of integers.

The following  theorem  is the main result  of this paper.
\begin{thm}.
 Let $f : \Rn\to\Co$  be the characteristic function of a probability measure $\mu$. Suppose that  $\mu$ has a nonzero absolutely continuous part  with continuous  density $\varphi$.   If  there exist $a\in\Rn$ and $\tau\in\Rn$, $\tau_k> 0$, $k=1,\dots,n$, such that
 \begin{equation}
 \tau\in (1/\sigma)Q_{2\pi}^n
\end{equation}
 and
 \begin{equation}
  (\su\, \varphi)\cap \bigl(a+\tau\,\Zn\bigr)=\emptyset,
\end{equation}
then $f$ has the $\sigma$-deterministic property.
 \end{thm}
The statement of this  theorem together with proposition 1.1   shows that the $\sigma$-deterministic   property of $f$ depends not  only on the size of both  $V_{\sigma}$ and  ${{\text{\,supp}}\,}\varphi$, but also on certain arithmetic properties of ${{\text{\,supp}}\,}\varphi$.

Now we will explain how accurate are our conditions (1.2) and  (1.3).   Given $\delta>0$ and $2\le q\le\infty$, let
\[
E^n_q(\delta)=\{x\in\Rn: \sum_{k=1}^n |x_k|^q<\delta\}.
\]
 Of course, $E^n_{\infty}(\delta)={\text{\rm{ int}}}Q^n_{(\delta,\dots,\delta)}$. Let $\varphi$ be any continuous density  such that
\begin{equation}
S(\varphi)=E^n_q(\delta).
\end{equation}
Take an arbitrary  $\sigma>0$ such that
\begin{equation}
0<\delta\sigma<\pi n^{\frac{1}{q}}.
\end{equation}
Then it is easily  checked that  there exist $a\in\Rn$ and $\tau\in\Rn$ such that (1.2) and (1.3) are satisfied.
More precisely, if (1.5) is satisfied, then there exists   $\varepsilon>0$ such that we can take $\tau_k=2\delta n^{-1/q}+\varepsilon$ for all $k=1,\dots, n$, and
\[
a=\delta n^{-\frac{1}{q}}(1+\varepsilon,\dots,1+\varepsilon).
\]
Hence,  (1.5) is  an example of  sufficient conditions for the $\sigma$-deterministic  property of $f$ with density $\varphi$ satisfying (1.4). It turns out that, at least for $q=2$,  this condition is sharp in the following sense.
\begin{thm}.
 Let $f : \Rn\to\Co$  be the characteristic function of a continuous probability density $\varphi$. Suppose that there exists $\delta>0$ such that
 \begin{equation}
S(\varphi)=E^n_2(\delta).
\end{equation}
 If
 \begin{equation}
\delta\sigma>\pi n^{\frac{1}{2}},
\end{equation}
$k=1,\dots,n$,   then  $f$  has not the $\sigma$-deterministic  property.
   \end{thm}
Let us remind of   our previous paper \cite{4}. The main result of \cite{4} states that if $\varphi$ is a continuous probability density of one variable and  $\Lambda_j=\tau_j+\alpha_j\Z$, $\tau_j\in\R$, $\alpha_j>0$,  $j=1,2$ are lattices with $\Lambda_1\cap \Lambda_2=\emptyset$ such that  $\varphi$ vanishes on  $\Lambda_1\cup \Lambda_2$ and $\alpha_j\sigma\le 2\pi$, $j=1,2$,  then, for any characteristic function $g:\R\to\Co$ with $g=\widehat{\varphi}$ on  $V _{\sigma}=\{x\in\R: |x|>\sigma\}$, we have that $g\equiv\widehat{\varphi}$, i.e.,  $\widehat{\varphi}$ has the $\sigma$-deterministic   property.

\section{ Preliminaries}\label{sec:2}
The aim of this section is to prove a number of auxiliary results on the theory of entire functions of several variables that are needed for the proofs of our main results.

For a closed subset $\Omega\subset\Rn$,  a function ƒ$\omega:\Rn\to \Co$     is called bandlimited to ƒ$\Omega$  if ƒ$\widehat{\omega}$  vanishes outside $\Omega$. The  subset $B^1_{\Omega}$ in  $L^1(\Rn)$ of  all $F\in L^1(\R)$ such that $F$ is bandlimited to  ƒ$\Omega$ is called the Bernstein space.     $B^1_{\Omega}$ is a Banach space with respect to the $ L^1(\Rn)$  norm. Below  we consider only the case where  $\Omega=\qu$. By the Paley-Wiener theorem,  any   $F\in \ben$ is infinitely differentiable on $\Rn$ and  has an extension onto the complex space $\Con$  to an entire function. For any  $F\in B^1_{\qu}$ and $\nu\in\Rn$, $\nu_k\neq 0$, $k=1,\dots, n$,  the Poisson summation formula (see, e.g., \cite[p. 166]{5})
 \begin{equation}
\sum_{\omega\in\nu\Zn}F(x+ \omega)=\frac1{\prod_{k=1}^{n}|\nu_k|}\sum_{\theta\in (1/\nu)\Zn}\widehat{F}\Bigl(2\pi \theta\Bigr)e^{2\pi i(x, \theta)}
\end{equation}
holds for all $x\in\Rn$. Note that both sums in (2.1) converge absolutely.

For $m=1,2,\dots,n$, let us denote by $H^n_m$ the set
\[
H^{n}_{m}=\{z=(z_1,\dots,z_n)\in\Con:  \ z_m\in\Z\}.
\]
It follows immediately that
\begin{equation}
k+ H^{n}_{m} = H^{n}_{m}
\end{equation}
 for any $k\in \Zn$ and if $\rho$ is a permutation of  $\{1,2,\dots, n\}$,  then
\begin{equation}
\lambda\in H^{n}_{m} \quad  {\text{\rm{if and only if }}}\quad (\lambda_{\rho(1)}, \lambda_{\rho(2)},\dots, \lambda_{\rho(n)})\in  H^{n}_{\rho(m)}.
\end{equation}
Set
\[
H^n=\bigcup_{m=1}^n H^n_m.
\]
In the sequel, we will use several times  the function
 \begin{equation}
\sn (z)=\prod_{k=1}^n \sin(\pi  z_k)
\end{equation}
with $z=(z_1,\dots,z_n)\in\Con$.  Obviously, $\sn $ is an  entire  function.  We denote by  $N(F)$  the zero set of an entire function $F$. Clearly,
\[
N(\sn)=H^n.
\]
Let $A\subset\Rn$. For  a function $g:A\to\R$, we say that $g$  has the same sign on $A$ if $g(t)\ge 0$ for all  $t\in A$, or $g(t)\le 0$ for all  $t\in A$.
\begin{prop}.
 Let $F\in\ben$. Assume that   there exist $a\in\Rn$ and $\tau\in\Rn$, $\tau_k>0$,  $k=1,\dots, n$,  such that
 \begin{equation}
F(x)=0 \quad {\text{\rm{and}}}\quad  \frac{\partial F(x)}{\partial x_k}=0
\end{equation}
for $k=1,\dots,n$,  and all $x\in a+\tau\Zn$.  Suppose, in addition,  that $F$ is real-valued on $\Rn$ and for each $x\in a+\tau\Zn$, there is  its neighborhood $M_x$ in $\Rn$ such that  $F$ has the same sign on $M_x$. If   (1.2) is satisfied, then $F\equiv 0$.
 \end{prop}
The proof of Proposition 2.1 is based on the following two lemmas.
\begin{lem}.
Let $F$ be an entire function on $\Con$  such that
\begin{equation}
F(z)=0 \quad {\text{ for \ all}}\quad z\in H^n_m, \ m=1,\dots, n.
\end{equation}
 Then there exists an entire function $G$ on $\Con$ such that
 \begin{equation}
F=G \cdot \sn.
\end{equation}
 \end{lem}

{\bf{Proof.}} \
 First, we claim that  there exists an entire function $f$ on $\Con$ such that
\begin{equation}
F(z)=\Bigl(\prod_{k=1}^n z_k \Bigr)f(z),
\end{equation}
$z\in\Con$. To that  end, we use the fact that any entire function $F$ can be expanded in a series in homogeneous polynomials
\begin{equation}
F(z)=\sum_{j=0}^{\infty} P_j(z),\quad P_j(z)=\sum_{|r|=j}c_rz^r,\quad z^r=\prod_{k=1}^nz_k^{r_k},
\end{equation}
where $r$ is a non-negative  $n$-multi-index,  $|r|=r_1+\dots + r_n$. Moreover, this series  converges uniformly on compact subsets of $\Con$. Set
\[
E^n_m=\{z\in\Con: \ z=(z_1,\dots,z_{m-1},0,z_{m+1},\dots,z_n)\},
\]
$m=1,\dots,n$.  Clearly,  $E^n_m \subset H^n_m$, $m=1,\dots,n$.  Hence, (2.6) is satisfied also for all $z\in E^n_m$. Next,  $E^m_n $  may be identified in a natural way with $\Co^{n-1}$.    Fix such an $m$. Then, using the  identity theorem for entire functions on $\Co^{n-1}$, we see that the condition $F=0$ on $E^n_m$, is equivalent to $c_r=0$ for each $r$ such that $r=(r_1,\dots,r_{m-1},0,r_{m+1},\dots,r_n)$. Therefore, $P_0=0$ and $P_j(z)=z_mQ_j(z)$, $j=1,\dots $, where  $Q_j$ is a homogeneous polynomial of degree $j-1$ or $Q_j\equiv 0$. By repeating the same for other $E^n_m$, we obtain that $P_1=\dots =P_{n-1}\equiv 0$ and
\begin{equation}
P_j(z)=\Bigl(\prod_{k=1}^n z_k \Bigr)R_j(z)
\end{equation}
for $j=n,n+1,\dots$, where $R_j$ is a homogeneous polynomial of degree $j-n$ or $R_j\equiv 0$. Next, the series $\sum_{j=n}^{\infty}R_j$ converges on compact subsets of $\Con$, since $\sum_{0}^{\infty}P_j$ has this property. Therefore, $\sum_{j=n}^{\infty}R_j$ defines  an entire function, say $f$,  on $\Con$. This, in light of (2.9) and (2.10),  proves our claim (2.8).

Now using the equality $H^n=N(\sn)$, we see that the function
\begin{equation}
G(z):=\frac{F(z)}{\Delta_n( z)}
\end{equation}
is well-defined for all $z\in \Con\setminus H^n$. We claim that $G$ can be  extended to an entire function on $\Con$. To this end, we  are going to use the Riemann removable singularity theorem (see, e.g. [6, p. 175]) in the case of this function  $G$ and the analytic set $H^n=N(\sn)$. Therefore, it is enough to show that $G$ is locally bounded on $H^n$, i.e., for every $\lambda\in  H^n$  there is its  open neighborhood  $ U_{\lambda}$ in $\Con$ such that  $G$ is bounded on $(\Con\setminus H^n)\cap U_{\lambda}$.

 Let $\lambda\in H^n$. Hence, $\lambda$ has at least one coordinate $\lambda_i$ such that $\lambda_i\in\Z$. Suppose that $\lambda$ has exactly $p$, $1\le p\le n$, such an $\lambda_{i_1},\dots,\lambda_{i_p}\in \Z$. Note that,   for any $\omega\in\Zn$ and for each perturbation $\rho$ of $\{1,2,\dots, n\}$,   we have
\[
\sn(z_{\rho(1)}+\omega_1,\dots, z_{\rho(n)}+\omega_n)=(-1)^{|\omega|}\sn( z)
\]
for all $z\in\Con$. Moreover, the function
\[
F_{\omega, \rho}(z):=F(z_{\rho(1)}+\omega_1,\dots, z_{\rho(n)}+\omega_n)
\]
also satisfies (2.6). Therefore, we can assume without loss of generality that
\begin{equation}
\lambda=(0,\dots,0, \lambda_{p+1},\dots, \lambda_n)
\end{equation}
with $\lambda_{p+1},\dots, \lambda_n\in \Co$, but  $\lambda_{p+1},\dots, \lambda_n\not\in\Z$.

Substituting (2.8) in (2.11),  we can rewrite (2.11) as
\begin{equation}
G(z)=\frac{f(z)}{a(z)b(z)}
\end{equation}
with
\begin{equation}
a(z)=\prod_{j=1}^{p}\frac{\sin(\pi z_j)}{z_j}\quad {\text{\rm{and}}} \quad b(z)=\prod_{j=p+1}^n\sin (\pi z_j).
\end{equation}
Of course, if $p=n$, then we take $b\equiv 1$. Clearly, $a$ and $b$ are entire functions. Moreover, using (2.12), we see that
\[
a(\lambda)=\pi^p \quad {\text{\rm{and}}} \quad b(\lambda) \neq 0.
\]
Hence, there exists $\varepsilon>0$ and a  neighbourhood $U_{\lambda}\subset \Con$ of $\lambda$ such that
\[
|a(z)|\ge \varepsilon \quad {\text{\rm{and}}} \quad |b(z)|\ge \varepsilon
\]
for all $z\in U_{\lambda}$. Now using the expression (2.13) and the fact that $f$ is entire, we obtain  that $G$ is bounded on $(\Con\setminus H^n)\cap U_{\lambda}$.
This yields our  claim. Then  (2.11) finishes the proof of Lemma 2.2.

\begin{lem}.
Let $G$ be an entire function on $\Con$ and suppose that  $G$ is real-valued on $\Rn$. Assume  that,  for each $x\in \Zn$, there is  its neighborhood $V_x$ in $\Rn$ such that  the function
\begin{equation}
F(z)=G(z)\sn( z)
\end{equation}
 has the same sign on $ V_x$. Then there exists an entire function $g$ such that
 \begin{equation}
G(z)=g(z)\sn( z)
\end{equation}
for all $z\in\Con$.
 \end{lem}
{\bf{Proof.}} \ We claim that
\begin{equation}
G(z)=0\quad {\text{for \ all}}\quad  z\in H^n_m, \ m=1,\dots,n.
\end{equation}
Fix any $x=(k_1,\dots,k_n)\in\Zn$. Take any number $\varepsilon$ such that
\begin{equation}
0<\varepsilon <1 \quad{\text{and}}\quad x+Q^n_{\varepsilon}\subset V_x,
\end{equation}
where $V_x$ is a neighborhood which satisfies  the hypothesis of our  lemma. Then, for any $t\in(-\varepsilon, \varepsilon)\setminus\{0\}$, we have that
\[
x_{\varepsilon,t}=(k_1+\varepsilon,\dots,k_{n-1}+\varepsilon, k_n+t)\in V_x
\]
 and $\sn(x_{\varepsilon,t})=-\sn(x_{\varepsilon,-t})$. Therefore, $G(x_{\varepsilon,t})=-G(x_{\varepsilon,-t})$, since the function (2.15)  has the same sign on $ V_x$. By the continuity of $G$,  we obtain that
  \begin{equation}
G(x_{\varepsilon,0})=0
\end{equation}
for all numbers $\varepsilon$  that satisfy  (2.18). Next, using the same argument, we obtain also that
\begin{equation}
G(x_{-\varepsilon,0})=0.
\end{equation}
Set $G_1(\lambda)=G(\lambda,k_n)$, $\lambda\in\Co^{n-1}$. Then $G_1$ is an entire function on $ \Co^{n-1}$. Moreover,  (2.19) and (2.20) show that $G_1$ vanishes on a real neighborhood of
$(k_1,\dots,k_{n-1})$, i.e., on the set
\[
(k_1,\dots,k_{n-1})+Q^{n-1}_{\varepsilon}.
\]
Therefore, the identity theorem for analytic functions ( see e.g., [6, p.21]) states that $G_1\equiv 0$ on $ \Co^{n-1}$. Hence $G=0$ on $H^n_n$. In the same way we obtain that $G=0$ on the other $H^n_m$, $m=1,\dots,n-1$. This proved our claim (2.17).

Finally, (2.17) shows that $G$ satisfies the hypotheses of Lemma 2.2. The proof is therefore complete.

{\bf{ Proof of Proposition 2.1. }} \  We start with the observation that it is sufficient to consider  the case if $a=0$, $\tau=(1,\dots, 1)$ and $F\in\bep$, i.e., that the conditions  (1.3) are satisfied for all $x\in\Zn$. Indeed, it is clear that $F\in \ben$ if and only if  $F_{a,\tau}(x):=F(a_1+\tau_1 x_1,\dots, a_n+\tau_nx_n)\in B^1_{\tau\cdot Q^n_{\sigma}}$. Also (1.2) implies that $B^1_{Q^n_{\tau\sigma}}\subset \bep$. Therefore, we can
 consider without loss of generality only functions  $F$ from the larger space $\bep$ and such that
\begin{equation}
F(x)=0 \quad {\text{\rm{and}}}\quad  \frac{\partial F}{\partial x_j}(x)=0
\end{equation}
for $j=1,\dots,n$ and all $x\in\Zn$.

Our proof is by induction on the dimension $n$ of $\Con$. If  $n=1$, then we make use of the fact that if $F$ is a function of one variable  such that $F\in B^1_{2\pi}:=B^1_{[-2\pi,2\pi]}$,  then   $F$ can be expanded using the Whittaker-Kotel'nikov-Shannon  sampling  formula with derivatives  (see \cite[p. 60]{2})
\[
F(z)=\sum_{n\in \Z}\Bigl(F(n)+F'(n)(z-n)\Bigr) \Biggl(\frac{\sin(\pi(z - n))}{\pi(z - n)}\Biggr)^2,
\]
$z\in\Co$. By (2.21), this formula implies that $F\equiv 0$.

Suppose that the proposition holds for dimension $n-1$.  Let $F\in \bep$ satisfy the hypothesis of our proposition.   Fix any $k\in\Z$ and set $F_1(\lambda):=F(k,\lambda)$, $\lambda \in \Co^{n-1}$. Then, for each  $\lambda\in\Z^{n-1}$, $F_1$ satisfies (2.21) and there is a   neighborhood $M_{\lambda}$ of $\lambda$ in $\R^{n-1}$ such that  $F_1$ has the same sign on $M_{\lambda}$. Using the induction hypothesis,  gives $F_1= 0$ on $\Co^{n-1}$. It follows that $F=0$ on $H^n_1$. In the same way we obtain that $F=0$ also on $H^n_{m}$, $m=2,\dots,n$. Therefore, by Lemma 2.2,  $F$ admits the representation  (2.7). Moreover, the functions  $F$ and $G$ in (2.7) satisfy the conditions of Lemma 2.3. Thus,
\begin{equation}
F=g\cdot {\sn^2},
\end{equation}
where $ g$ is an entire function on $\Con$.

We claim that $F(z)=0$ for all $z\in \Con\setminus N(\sn)=\Con\setminus H^n$. Fix any $z\in \Con\setminus N(\sn)$. For $\lambda\in \Co$, let us define $F_z(\lambda)=F(\lambda,z_2,\dots,z_n)$. By (2.22), we get
\[
F_z(\lambda)=\sin^2\pi \lambda\cdot g(\lambda, z_2,\dots,z_n)\cdot\prod_{j=2}^n \sin^2\pi z_j .
\]
Then $F_z\in B^1_{2\pi}$, since  $F\in \bep$. Moreover,  $F_z$ satisfies (2.21) in the case of one variable. Then $F_z(\lambda)=0$ for all $\lambda\in\Co$, as was proved above. Therefore, $F(z)=F_z(z_1)=0$.  This yields the claim. 	

Finally, since $N(\sn)$ is a closed and proper subset of $\Con$, it follows that the entire function $F$ vanishes on an open and nonempty subset $\Con\setminus N(\sn)$ of $\Con$.  Thus, $F\equiv 0$ and Proposition 2.1 is proved.

\section{ Proofs}\label{sec:2}

We  define the inverse Fourier transform
 \[
\check{\chi}(t)=\frac1{(2\pi)^n}\int_{\Rn}e^{i(t, x)}\chi(x)\,dx
\]
so that   the inversion  formula  $\widehat{(\check{\chi})}=\chi$ holds for suitable $\chi\in L^1(\Rn)$.  Let $\Bn=\{\widehat{\mu}: \mu\in\Mn\}$ be the Fourier-Stieltjes algebra with the usual pointwise multiplication. The norm in  $\Bn$   is inherited  from  $\Mn$ in such a way $\|\hat{\mu}\|_{\Bn}:=\|\mu\|_{\Mn}$. The Fourier algebra $A(\Rn)=\{\widehat{\varphi}: \varphi\in L^1(\Rn)\}$  is an ideal in $\Bn$.

 As usual, we write $S(\Rn)$ for the  Schwartz space of test functions on $\Rn$ and $S'(\Rn)$ for the dual space of tempered distributions.

 \begin{prop}
Let  $ \mu,\eta\in\Mn$. Assume  that $\mu=\mu_a+\mu_s$ and $\eta=\eta_a+\eta_s$ are the usual Lebesgue decompositions of $ \mu$ and $\eta$ into their  absolutely continuous and singular parts, respectively. If there exists  an $V_{\sigma}$  such that $\widehat{\mu}=\widehat{\eta}$  on $V_{\sigma}$, then $\mu_s=\eta_s$.
\end{prop}
{\bf Proof.}  \ Let $u\in S(\Rn)$ be such that $u(t)=1$ for all $t\in Q_{\sigma}^n$.  Then
\begin{equation}
\widehat{\mu}-\widehat{\eta}\equiv u(\widehat{\mu}-\widehat{\eta}).
\end{equation}
Since $S(\Rn)\subset A(\Rn)$ and $A(\Rn)$ is an ideal in $B(\Rn)$, it follows from (3.1) that there exists $\chi\in L^1(\Rn)$ such that $\widehat{\mu}-\widehat{\eta}=\widehat{\chi}$. Therefore,  $\mu_s=\eta_s$ and Proposition 3.1 is proved.

{\bf{ Proof  of Theorem 1.2.}}  \  Assume  that  $g$ is any characteristic function such that $g=f$ on $\us$. We need to show that $f\equiv g$. By Proposition 3.1, without loss of generality, we may assume that $f=\widehat{\varphi}$. Then  we first claim that if $\mu_g$  is the probability measure such that $g=\widehat{\mu_g}$, then  $\mu_g$ is absolutely continuous with respect to the Lebesgue measure on $\Rn$. Indeed, using that $f-g$ is   compactly  supported on  $\Rn$ and $f-g\in B(\Rn)$, we have in  a similar way as in the proof of Proposition 3.1, that  $f-g\in A(\Rn)$. Since $f\in  A(\Rn)$, it follows that $g\in  A(\Rn)$, i.e., $\mu_g$  has a  density function $\vartheta$. Therefore, our claim is   proved.  Note that $\vartheta$ is also continuous. More precisely, $\vartheta\in  C_0(\Rn)$, where   $C_0(\Rn)$ is the usual space of continuous functions on $\Rn$ that vanish at infinity.  This follows from facts that both $\varphi$ and $\varphi-\vartheta$ are elements of $L^1(\Rn)$.

  Let us define
   \begin{equation}
\xi= \varphi-\vartheta.
 \end{equation}
Then  $\xi\in\bev$, since $g-f=0$ on $V_{\sigma}$.  Furthermore, (3.2) implies that
 \begin{equation}
\int_{\Rn}\xi(x)\,dx=0
 \end{equation}
 and
 \begin{equation}
\xi(x)\le  \varphi(x)
 \end{equation}
for all $x\in\Rn$.  Combining (1.2) with (3.4), we see that
 \begin{equation}
\xi(x)\le 0
 \end{equation}
for $x\in a+\tau\Zn$.

Next, we claim that
\begin{equation}
\xi(x)= 0,
 \end{equation}
for all $x\in a+\tau\Zn$. Indeed, applying (2.1) in the case of $F=\xi$, $x=a$, and $\nu=\tau$, we get
\begin{equation}
\sum_{\omega\in\tau\Zn}\xi(a+\omega)=\frac1{\prod_{k=1}^{n}\tau_k}\sum_{\theta\in(1/\tau)\Zn}\widehat{\xi}\Bigl(2\pi \theta\Bigr)e^{2\pi i(a, \theta)}.
\end{equation}
If $\theta\in \Zn$ and $\theta\neq 0$, then (1.2) implies that
\[
(2\pi)/\tau\in \partial \qu\quad {\text{\rm{ or}}}\quad (2\pi)/\tau \not\in \qu.
\]
Hence,
\[
\widehat{\xi}\bigl( (2\pi)/\tau\bigr)=0
\]
for $\theta\in\Zn$, $\theta\neq 0$,  since $\xi\in\bev$.  On the other hand,  (3.3) gives  that $\widehat{\xi}(0)=0$. Therefore,  we can rewrite (3.7) as
\[
\sum_{\omega\in\Zn}\xi(a+\tau \omega)=0.
\]
Combining this with (3.5), we obtain our claim (3.6).

Now,  if $a$ and $\tau$  are the same as in (1.3), then, for  each   $x=a+\tau\omega\in\Rn $, $\omega\in\Zn$,  there exists a neighborhood $M_x$ of $x$ in $\Rn$ such that $M_x\cap{\text{supp}} \,\varphi=\emptyset$. Therefore, (3.4) implies that $\xi\le 0$ on $M_x$. In addition, by (3.6), we see that $x$ is a local maximum point of $\xi$.  Hence, for an infinitely differentiable and real-valued on $\Rn$ function $\xi$, we get that
\[
\frac{\partial \xi(x)}{\partial x_k}=0
\]
for $k=1,\dots, n$ and all  $x\in a+\tau\Zn$. Using this and (3.6), we obtain by  Proposition 2.1,  that $\xi\equiv 0$. This is equivalent to $g\equiv f$. Theorem 1.2  is proved.

{\bf { Proof  of Theorem 1.3. }} \ The argument at the  beginning of the proof of Theorem 1.2 shows that it is sufficient to prove
that there exists $\xi\in\bev$, $\xi\not\equiv 0$, satisfying (3.3) and (3.4).  We claim that there exist  $\varrho>0$ and $A>0$ such that the function
\begin{equation}
\xi(x)=\varrho\prod_{k=1}^n\biggl(\frac{\cos \frac{\sigma x_k}{2}}{x_k^2-(\pi/\sigma)^2}\biggr)^2\Bigl[A^2-\sum_{k=1}^nx_k^2\Bigr]
\end{equation}
is an element of $\bev$ and satisfies (3.3) and (3.4). Set
\[
\xi_1(x)=\biggl(\frac{\cos \frac{\sigma x}{2}}{x^2-(\pi/\sigma)^2}\biggr)^2 \quad {\text{\rm{and}}} \quad \xi_2(x)= x^2\xi_1(x),
\]
$x\in\R$. It is obvious that   $\xi_1, \xi_2\in B^1_{Q^1_{\sigma}}$. Next,  (3.8) admits the following expansion
\begin{equation}
\xi(x)=\varrho\biggl( A^2\prod_{k=1}^n\xi_1(x_k) - \sum_{k=1}^n\Bigl[\xi_2(x_k)\cdot \prod_{j=1, \ j\neq k}^n \xi_1(x_j)\Bigr] \biggr).
\end{equation}
This implies that  $\xi\in \bev$.

Now we find such an $A$ for which (3.3) is satisfied. For this reason,   we apply (2.1)  to $\xi_m$, $m=1,2$ in the case of $n=1$, $x=\pi/\sigma$, and $\nu=2\pi/\sigma$. Then
 \begin{equation}
\sum_{l\in \Z}\xi_m\biggl(\frac{\pi}{\sigma}+\frac{2\pi l}{\sigma}\biggr)=\frac{\sigma}{2\pi}\sum_{k\in\Z}\widehat{\xi_m}(\sigma k)e^{i k\pi }.
\end{equation}
Since $\xi_m\biggl(\frac{\pi}{\sigma}+\frac{2\pi l}{\sigma}\biggr)=0$ for any $l\in\Z\setminus\{-1,0\}$  and $\widehat{\xi_m}(\sigma k)=0$ for each $k\in\Z\setminus\{0\}$,  it follows from (3.10) that
 \begin{equation}
\int_{\R}\xi_m(t)\,dt= \widehat{\xi_m}(0)= \frac{2\pi}{\sigma}\biggl(\xi_m\biggl(\frac{\pi}{\sigma}\biggr)+\xi_m\biggl(-\frac{\pi }{\sigma}\biggr)\biggr)=\frac{4\pi}{\sigma}\xi_m \biggl(\frac{\pi}{\sigma}\biggr).
\end{equation}
A straightforward calculation gives
 \[
\xi_1 \biggl(\frac{\pi}{\sigma}\biggr)= \frac{\sigma^4}{16\pi^2} \quad {\text{\rm{and}}} \quad \xi_2 \biggl(\frac{\pi}{\sigma}\biggr)=\frac{\pi^2}{\sigma^2}\xi_1 \biggl(\frac{\pi}{\sigma}\biggr)=\frac{\sigma^2}{16}.
\]
Hence, combining this with (3.9) and (3.11), we obtain
\begin{gather}
\int_{\Rn}\xi(x)\,dx= \varrho\biggl( A^2\prod_{k=1}^n\int_{\R}\xi_1(x_k)\,dx_k - \sum_{k=1}^n\Bigl[\int_{\R}\xi_2(x_k)\,dx_k\cdot \prod_{j=1, \ j\neq k}^n \int_{\R}\xi_1(x_j)\,dx_j\Bigr] \biggr)\nonumber\\
=  \varrho\biggl( A^2\biggl(\frac{\sigma^3}{4\pi}\biggr)^n  - n\frac{\sigma\pi}{4}\biggl(\frac{\sigma^3}{4\pi}\biggr)^{n-1}\biggr)=   \varrho\biggl( \frac{\sigma^3}{4\pi}\biggr)^n\biggl(A^2- \frac{\pi^{2}}{\sigma^2}n\biggr).   \nonumber
\end{gather}
Therefore, if we take
\[
A=\frac{\pi \sqrt{n}}{\sigma},
\]
then (3.3) is satisfied.

Finally, we show that there exists $\varrho>0$ such that the function (3.8) satisfies (3.4). Indeed, (3.4) is satisfied on
\[
\Rn\setminus E^n_2(A)=\Rn\setminus E^n_2\Bigl(\frac{\pi \sqrt{n}}{\sigma}\Bigr),
\]
since $\varphi\ge 0$ on $\Rn$ and $\xi\le 0$ on $\Rn\setminus E^n_2((\pi\sqrt{n})/\sigma)$. On the other hand, $\varphi$ is continuous on $\Rn$. Hence,  (1.6) and (1.7) imply that
\[
\inf_{x\in E^n_2((\pi\sqrt{n})/\sigma)} \varphi(x)>0.
\]
Using this and the fact that $\xi$ is a bounded function on $\Rn$, we see that there exists $\varrho>0$ such that  (3.4) is satisfied also on $E^n_2((\pi\sqrt{n})/\sigma)$. This finishes our proof.


}}
\end{document}